\documentclass[a4paper]{article}

\usepackage{graphicx}
                     
\usepackage{amsmath,amssymb}

\usepackage{stmaryrd}
\usepackage{booktabs}
\usepackage{a4wide}
\usepackage{fix-cm}
\usepackage{microtype}

\begin{document}

\title{A Model Reduction Framework for Efficient Simulation of Li-Ion Batteries}
\author{Mario Ohlberger\thanks{Center for Nonlinear Science \& Applied Mathematics Muenster, Einsteinstrasse 62, 48149 Muenster, Germany, \texttt{\{mario.ohlberger, stephan.rave\}@uni-muenster.de}} \and Stephan Rave\footnotemark[1] \and Sebastian Schmidt\thanks{Fraunhofer Institute for Industrial Mathematics ITWM, Fraunhofer-Platz 1, 67663 Kaiserslautern, Germany, \texttt{sebastian.schmidt@itwm.fraunhofer.de}} \and Shiquan Zhang\thanks{School of Mathematics, Sichuan University, Chengdu 610064, China, \texttt{shiquanz3@gmail.com}}}

\maketitle

\begin{abstract}
In order to achieve a better understanding of degradation processes in lithium-ion batteries, the
modelling of cell dynamics at the mircometer scale is an important focus of current mathematical
research. 
These models lead to large-dimensional, highly nonlinear finite volume discretizations which, due
to their complexity, cannot be solved at cell scale on current hardware.
Model order reduction strategies are therefore necessary to reduce the computational complexity
while retaining the features of the model.
The application of such strategies to specialized high performance solvers asks
for new software designs allowing flexible control of the solvers by the reduction algorithms.
In this contribution we discuss the reduction of microscale battery models with the reduced basis method
and report on our new software approach on integrating the model order reduction software pyMOR with third-party solvers. 
Finally, we present numerical results for the reduction of a 3D microscale battery model with porous
electrode geometry.
\end{abstract}

\section{Introduction}
\label{sec:introduction}

A major cause for the failure of rechargeable lithium-ion batteries is the deposition of metallic lithium at the
negative battery electrode (Li-plating).
Once established, this metallic phase can grow in the form of dendrites to the positive electrode,
ultimately short-circuiting the cell.
As Li-plating is initiated at the interface between active electrode particles and the electrolyte,
understanding of this phenomenon is only gained through physical models accounting for effects on
the micrometer-scale. This in turn requires highly resolved meshes in the model discretization.

A thermodynamically consistent microscale battery model was developed in \cite{LatzZausch2011}.
Based on a finite volume discretization \cite{PopovVutovEtAl2011}, this model has been implemented
at Fraunhofer ITWM in the battery simulation software BEST \cite{LessSeoEtAl2012}.  
However, since such microscale discretizations lead to very large, highly nonlinear equation systems,
simulations can currently only be performed on small portions of the cell and parameter studies
testing different charging regimes or operating conditions are very time consuming.
It is therefore desirable to combine microscale modeling with model order reduction
strategies which are able to reduce the computation time while at the same time keeping the
microscopic features of the model.

The reduced basis method is a well-established approach for model order reduction of problems given 
by parametric partial differential equations and has been successfully adapted to various industrial
applications (see references in \cite{HaasdonkOhlberger2008}).
In this approach, the original equation is projected onto a low-dimensional discrete function space which has been constructed from the solution
trajectories of the high-dimensional problem for selected parameters of a well-chosen training set.
The applicability of the method to nonlinear finite volume discretizations has been been shown in
\cite{HaasdonkOhlberger2008,DrohmannHaasdonkEtAl2012}.
Results for the model order reduction of a pseudo-2D battery model using similar techniques have been presented in \cite{IlievLatzEtAl2012}.

A major challenge for the implementation of reduced basis schemes lies, however, in their integration with
(already existing) PDE solvers:
in those schemes the solver has to be controlled by the reduction algorithm which, apart from
solving the high-dimensional problem, now also has to provide the reduction data needed to
perform the low-dimensional simulations. 
Moreover,
the solver is usually unable to perform the reduced computations, which are based on different data
structures.
This often leads to insertion of model reduction specific algorithms into the solver's code base, while in a separate code base the
solution algorithm for the reduced problem is re-implemented \cite{DrohmannHaasdonkEtAl2012a}.
As a result, code is duplicated and the adoption of a different model reduction strategy requires
changes in both code bases.

After discussing the application of the reduced basis method to the microscale model from \cite{LatzZausch2011}, we
present the design of our new model reduction software {\nobreak pyMOR} \cite{PYMOR} which is
specifically tailored to address these problems by offering a deep and flexible integration with
external PDE solvers.
We will conclude with first numerical results for the reduction of the full 3D-model with porous electrode geometries,
underlining the potential of the model reduction approach.

\section{Reduction of the Microscale Model}
\label{sec:model}
Our work is based on the microscale battery model introduced in \cite{LatzZausch2011}.
Under the assumption of a globally constant temperature $T$, this model is given by a system of
partial differential equations for the concentration of Li$^+$-ions $c$ and the electrical potential
$\phi$ on each part of the domain, i.e. the positive and negative electrodes, the
electrolyte and the current collectors.
Each of these systems is of the form
\begin{equation*}
        \frac{\partial c}{\partial t} + \nabla \cdot N = 0, \qquad 
        \nabla \cdot j = 0,
\end{equation*}
where $N = -(\alpha(c, \phi) \nabla c + \beta(c, \phi)\nabla \phi)$, $j = -(\gamma(c, \phi) \nabla c +
\delta(c, \phi) \nabla \phi)$
with the coefficients $\alpha, \beta, \gamma, \delta$ depending on the domain for which the system
is given.
While these coefficients can be considered constant in first approximation, a strong nonlinearity
enters the model through the interface conditions between electrolyte and active particles in the
electrodes.
These conditions are given by prescribing the normal interface fluxes of concentration and potential
into the electrolyte via the Butler-Volmer kinetics, i.e. 
\begin{eqnarray*}
       j_s \cdot n &=& j_e \cdot n = 2k\sqrt{c_ec_s(c_{max}- c_s)}
                \sinh\left(\frac{\phi_s - \phi_e - U_0(\frac{c_s}{c_{max}})}{2RT} \cdot F \right),\\
\end{eqnarray*}
and $N_s \cdot n = N_e \cdot n = j_{s} \cdot n / F$.
Here the subscripts $s$ ($e$) denote the value of the respective quantity in the active particle
(electrolyte) domain at the interface, and $n$ is the unit normal at the interface pointing into the
electrolyte. $U_0$ denotes the open circuit potential, $k$ is a reaction rate,
$c_{max}$ the maximum Li-ion concentration in the particle and
$T$ the temperature. The constants $F$ and $R$ denote the Faraday and universal gas constants.
The system is closed via appropriate boundary conditions as well as interface conditions for the current collectors.
E.g. a constant charge rate $I$ corresponds to the Neumann boundary condition $j \cdot n = -I$ at
the positive electrode side of the domain.

\begin{figure}[t]
\begin{center}
	\includegraphics[scale=0.17, trim=0 80 0 0, clip]{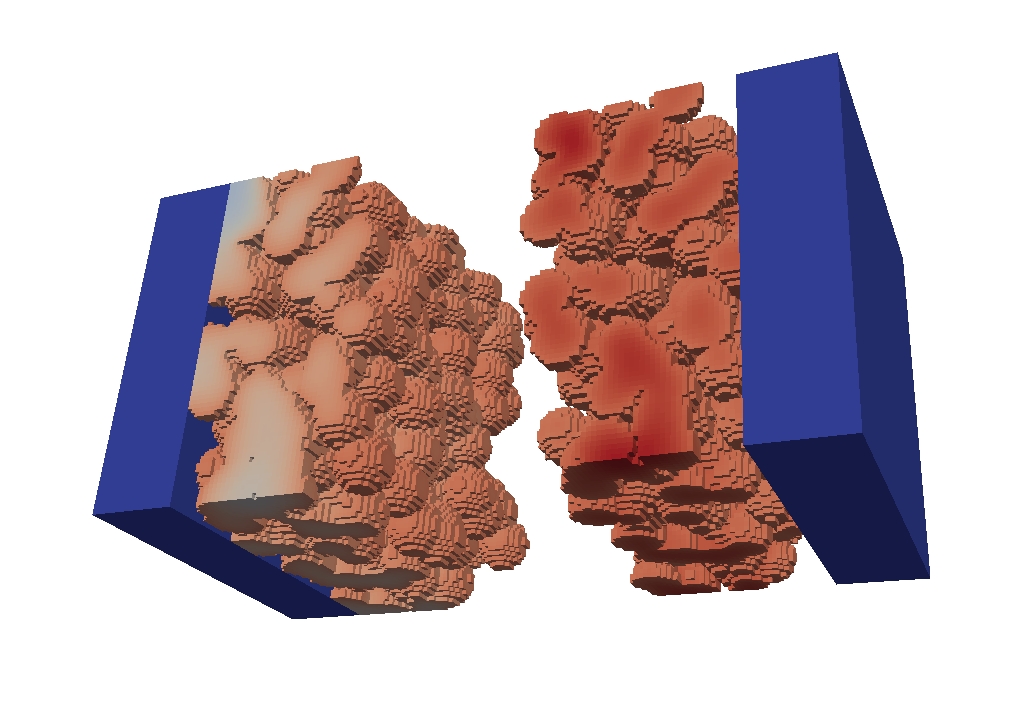}\\
\end{center}
        \caption{%
                Detailed simulation of battery model with DUNE on a $48\mu m
		\times 24\mu m \times 24\mu m$ computational domain with random electrode geometry.
		Coloring indicates Li$^+$ concentration in active particles (electrolyte not
		displayed).}
	\label{fig:simulation}
\end{figure}
\subsection{Discretization}
A discretization of the model based on a cell centered finite volume scheme has been introduced in
\cite{PopovVutovEtAl2011}. In this discretization, the interface conditions between electrolyte and
active particles are incorporated into the numerical fluxes and the implicit Euler method is used
for time discretization. As a result, one obtains nonlinear equation systems of the form

\begin{equation}
\label{eq:detailed}
        \begin{bmatrix}
                 \frac{1}{\Delta t}(c_{\mu}^{(t+1)} - c_{\mu}^{(t)}) \\
                 0
        \end{bmatrix}
         + A_\mu
         \left(\begin{bmatrix}
                c_{\mu}^{(t+1)} \\
                \phi_{\mu}^{(t+1)}
         \end{bmatrix}\right)
         = 0, \qquad c_{\mu}^{(t)}, \phi_{\mu}^{(t)} \in V_h
\end{equation}
with $A_\mu$ denoting the finite volume space operator acting on the discrete function space
$V_h \oplus V_h$.
The subscript indicates the dependence of the solution on a certain set of parameters $\mu$
(we consider the charge rate and temperature in our example below).
The discrete equation systems are solved in BEST with a Newton scheme utilizing an algebraic
multigrid solver for the linear systems in each Newton step.

\subsection{Reduced Basis Approximation}
\label{sec:rb}
The reduced basis method is based on the idea of performing a Galerkin projection of the
high-dimensional discrete equations (\ref{eq:detailed}) onto low-dimensional
subspaces $\tilde{V}_c, \tilde{V}_\phi \subset V_h$ constructed from solutions of (\ref{eq:detailed}) for
appropriately selected parameters.
Under this projection, (\ref{eq:detailed}) is transformed into
\begin{equation}
\label{eq:reduced}
                \begin{bmatrix}
                         \frac{1}{\Delta t}(\tilde{c}_{\mu}^{(t+1)} - \tilde{c}_{\mu}^{(t)}) \\
                         0
                \end{bmatrix}
                 + 
                \left\{
                 P_{\tilde{V}} \circ A_\mu
         \right\}
                 \left(\begin{bmatrix}
		        \tilde{c}_{\mu}^{(t+1)} \\
			\tilde{\phi}_{\mu}^{(t+1)}
                 \end{bmatrix}\right)
         = 0, \qquad \tilde{c}_{\mu}^{(t)} \in \tilde{V}_c,\, \tilde{\phi}_{\mu}^{(t)} \in \tilde{V}_\phi,
\end{equation}
where $P_{\tilde{V}}$ denotes the orthogonal projection onto the reduced space $\tilde{V}:=
\tilde{V}_c \oplus \tilde{V}_\phi$.
After this projection has been performed in a preceding ``offline-phase'', the resulting
low-dimensional system can be solved quickly for new parameter values in a following 
``online-phase''.

For the selection of $\tilde{V}_c$ and $\tilde{V}_\phi$ a large variety of algorithms
has been considered (\cite{HaasdonkOhlberger2008} and references therein), many of which are based
on a greedy search over a prescribed (or adaptively refined) training set of parameters:
in each round of the algorithm, an error estimator is used to search the training set for the
parameter $\mu^*$ to which the solution of (\ref{eq:detailed}) is worst
approximated by the solution of the reduced problem (\ref{eq:reduced}).
The high-dimensional solution trajectory $[c^{(t)}_{\mu^*}, \phi^{(t)}_{\mu^*}]$ is then computed and
$\tilde{V}_c$, $\tilde{V}_\phi$ are enlarged by vectors from the linear span
of this trajectory via an appropriate extension algorithm.
As the reduced spaces are constructed from solutions of the full microscale model, characteristic
features, e.g. concentration hotspots in certain electrode regions due to local particle
geometry, are still representable within these spaces, despite their low dimensionality.

While posed on low-dimensional spaces, problem (\ref{eq:reduced}) still depends on evaluations of the
high-dimensional operator $A_\mu$.
This dependency can be removed by application of the so-called empirical operator interpolation
method \cite{DrohmannHaasdonkEtAl2012}.
In this approach, the given operator is only evaluated at a small number of degrees of freedom (DOFs)
of the discrete space.
The evaluation of the full operator is then approximated via linear combination with a
pre-computed (collateral) interpolation basis.
The interpolated operator can be evaluated quickly, independently of the dimension of $V_h$,
due to the locality of finite volume operators:
the evaluation of $A_\mu$ at $M$ degrees of freedom only requires the knowledge of its argument at $M^\prime \leq 
C \cdot M$ DOFs with $C$ being determined by the maximum number of cell neighbours in the given grid.
If we denote by $\tilde{A}_\mu: \mathbb{R}^{M^\prime} \to \mathbb{R}^{M}$ the restricted operator
and by $R_{M^\prime}: V_h^2 \to \mathbb{R}^{M^\prime}$, $I_M: \mathbb{R}^{M} \to
V_h^2$ the operators given by projection onto the interpolation DOFs and linear combination
with the collateral basis, we obtain the fully reduced equation systems

\begin{equation}
\label{eq:interpolated}
                \begin{bmatrix}
                         \frac{1}{\Delta t}(\tilde{c}_{\mu}^{(t+1)} - \tilde{c}_{\mu}^{(t)}) \\
                         0
                \end{bmatrix}
                 + 
                \left\{
                (P_{\tilde{V}} \circ I_{M}) \circ
                \tilde{A}_\mu \circ R_{M^\prime}
         \right\}
                 \left(\begin{bmatrix}
			\tilde{c}_{\mu}^{(t+1)} \\
			\tilde{\phi}_{\mu}^{(t+1)}
                 \end{bmatrix}\right)
         = 0.
\end{equation}

The linear operators $P_{\tilde{V}} \circ I_M$ and $R_{M^\prime}$ can be pre-evaluated during the
offline-phase for a given basis of $\tilde{V}$, completely eliminating high-dimensional operations from
(\ref{eq:interpolated}).
For the determination of the interpolation DOFs and collateral basis, greedy search strategies can
again be utilized \cite{DrohmannHaasdonkEtAl2012}.

\section{A new Software Framework}
\label{sec:pyMOR}

The implementation of reduced basis schemes involves several building blocks: solution of the
detailed problem (\ref{eq:detailed}) for a given parameter, projection of the operators,
extension of the reduced spaces (high-dimensional operations), as well as solution of the reduced
problem (\ref{eq:interpolated}), estimation of the reduction error and greedy algorithms
(low-dimensional operations).
In previous software approaches \cite{DrohmannHaasdonkEtAl2012a}, the implementation of all
high-dimensional operations takes place in the solver code, whereas the low-dimensional operations are
implemented in a separate model reduction software.
As a consequence, 
both code bases have to be adapted if the reduction strategy shall be modified.
This can slow down implementation of new algorithms significantly if the solver is developed by a different team than the model reduction software.
Moreover, despite the fact that
(\ref{eq:detailed}) and (\ref{eq:interpolated}) are of the same mathematical structure, both software packages need to implement the same algorithm for solving the respective problems.
In particular, for empirical operator interpolation the restricted operator $\tilde{A}_\mu$ has to be implemented again for the reduced
scheme.

The design of pyMOR mitigates these difficulties by exploiting the observation that all
aforementioned building blocks can be implemented in terms of operations on the following types of objects, either provided by implementations in {\nobreak pyMOR} itself (usually low-dimensional objects) or by external solvers (usually high-dimensional objects):
\begin{itemize}
\item \textbf{Vector arrays} store collections of vectors, supporting basic linear algebra operations, e.g.
	computation of linear combinations of vectors or scalar products.
	Selected DOFs can be extracted for the implementation of operator interpolation.
\item \textbf{Operators} represent linear or nonlinear operators, bilinear forms or functionals.
        Operators can be applied to vector arrays. Linear solvers are exposed through application
	of the inverse operator, Jacobians and restricted operators can be formed. 
\item \textbf{Discretizations} encode as containers for operators the mathematical structure of a
        given discrete problem and implement algorithms for solving the problem in terms of the operators they contain.
\end{itemize}
All algorithms in pyMOR are implemented in terms of the interfaces provided by these classes. As an
important consequence, there is no distinction between high- and low-dimensional objects in pyMOR except for
the different types of vector arrays or operators that represent them. In particular, the same
discretization class
can be used to solve (\ref{eq:detailed})
as well as (\ref{eq:interpolated}) or (\ref{eq:reduced}). The reduction process
merely consists in the replacement of
operators of a given discretization object by the corresponding projected operators.
For empirical interpolation, pyMOR implements a generic interpolated operator which can be used to efficiently interpolate any restrictable operator in pyMOR.
The evaluation of the restricted operator $\tilde{A}_\mu$ can still be performed by the same code used to evaluate the full operator $A_\mu$.

\begin{figure}[t]
\begin{center}
\includegraphics{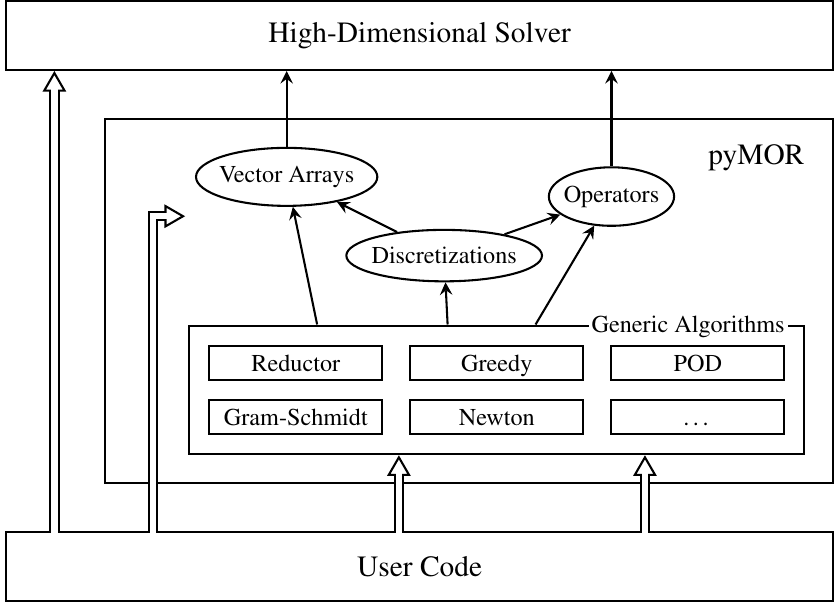}
\caption{Sketch of the interface concept for the integration of pyMOR with external solvers.}
\label{fig:design}
\end{center}
\end{figure}

As a consequence of this design, the model reduction algorithms in pyMOR are completely decoupled from
the development of the high-dimensional discretizations
(cf. Fig.~\ref{fig:design}).

\subsection{Implementational aspects}
Following the line of most other model order reduction packages, we chose with Python
a scripting language for the implementation of
pyMOR. Such languages offer a high amount of interactivity, making it very easy to experiment with 
various variants of model reduction algorithms.

While there is no underlying assumption of how the communication through the abstract interfaces is handled,
we favour, where possible, a tight integration of external solvers with pyMOR.
In particular for shared-memory solvers, an attractive option is the compilation of the
solver code as a shared library which then can be directly loaded as a Python extension module.
Apart from offering the easiest and at the same time most efficient way of
integration, an additional benefit is the direct accessibility of solver data structures from Python
which can be exploited to quickly augment the high-dimensional code with additional features.
This route of development has also been chosen for the ongoing integration of pyMOR with BEST within
the publicly founded MULTIBAT project.

\section{Numerical Results}

\begin{table}[t]
\caption{Constants used in numerical example, $c_0$ denotes initial concentration. Furthermore, $U_0(x) =
	-0.132 + 1.41\cdot \exp(-3.52x)$ for the negative and $U_0(x) = 4 + 0.07\cdot
	\tanh(-22x + 12) - 0.1\cdot(1/(1.002 - x)^{0.37} - 1.6) - 0.045\cdot \exp(-72 x^8) +
	0.01\cdot \exp(-200(x - 0.19))$ for the positive electrode, $R = 8.314$, $F =
	9.6487 \cdot 10^4$.}
\label{tab:constants}
\setlength{\tabcolsep}{1.7mm}
\begin{center}
\begin{tabular}{lccccccc}
\toprule
domain & $\alpha$ & $\beta$ & $\gamma$ & $\delta$ & $c_0$ & $c_{max}$ & k \\
\midrule
electrolyte    & $1.622\cdot 10^{-6}$   & $0$ &  $-5.171\cdot 10^{-5} \cdot T$ & $0.02$ & $1.200 \cdot 10^{-3}$   & -- & -- \\
pos. electrode & $1.0 \cdot 10^{-10}$   & $0$ &  $0$                           & $0.38$ & $2.057 \cdot 10^{-2}$   & $2.367\cdot 10^{-2}$ & 0.2 \\
\ \  --- current coll. & $0$                    & $0$ &  $0$                           & $0.38$ & $0$                     & -- & -- \\
neg. electrode & $1.0 \cdot 10^{-10}$   & $0$ &  $0$                           & $10$   & $2.639 \cdot 10^{-3}$   & $2.468\cdot 10^{-2}$ & 0.002\\
\ \  --- current coll. & $0$                    & $0$ &  $0$                           & $10$   & $0$                     & -- & --\\
\bottomrule
\end{tabular}
\end{center}
\end{table}

In order to provide a testbed for our reduction framework, an experimental implementation of the battery model has been
developed based on the PDELab discretization module for the DUNE software
framework \cite{BastianBlattEtAl2008} (cf. Fig.~\ref{fig:simulation}).
As a first experiment, we considered a small 3D test problem with randomly generated
electrode geometry, for which we evaluated the approximation quality of the reduced basis projection
(\ref{eq:reduced}).
We chose constant material properties resulting in the coefficients in Table \ref{tab:constants}.
The computational domain was of size $4.8\cdot 10^{-3} \times 2.4\cdot 10^{-2} \times 2.4\cdot 10^{-2}$
($cm^3$) which was meshed with a regular $40\times 20 \times 20$ grid. The width of the electrodes
(current collectors) was 10 (5) grid cells. The positive (negative) electrode was filled to 61.4\%
(74.2\%) with particle cells. 20 time steps of length 30 (s) were made.

\begin{table}[t]
\caption{Relative $L^\infty-L^2$ errors for the reduced basis approximation (\ref{eq:reduced}) of the high-dimensional model (\ref{eq:detailed}).
The basis size denotes the (equal) dimensions of $\tilde{V}_c$ and $\tilde{V}_\phi$.}
\label{tab:errors}
\setlength{\tabcolsep}{4.5mm}
\begin{center}
\begin{tabular}{lcccc}
\toprule
basis size & 8 & 16 & 24 & 32 \\
\midrule
concentration & $8.7 \cdot 10^{-3}$ & $1.9 \cdot 10^{-3}$ & $1.2 \cdot 10^{-3}$ & $4.3 \cdot 10^{-4}$ \\
potential & $1.3 \cdot 10^{-3}$ &  $2.1 \cdot 10^{-4}$ &  $7.7 \cdot 10^{-5}$ &  $1.5 \cdot 10^{-5}$ \\
\bottomrule
\end{tabular}
\end{center}
\end{table}

The parameters, charge rate $I$ and temperature $T$, were allowed to vary in
the intervals $[10^{-4},\ 10^{-3}]$ ($A/cm^2$) and $[250,\ 350]$ ($K$).
The reduced spaces were constructed with the POD-Greedy algorithm \cite{HaasdonkOhlberger2008} on a training set of $3 \times 3$
equidistant parameters, using the true reduction error for snapshot selection. During each extension
step, both reduced spaces were extended separately by orthogonally projecting the selected trajectory onto the respective reduced space and then enlarging the space with the
first POD mode of the trajectory of projection errors.
In Table \ref{tab:errors}, the maximum reduction error over the whole parameter space is estimated for different basis sizes 
by computation of the errors for 20 randomly selected new parameters.

\subsection*{Acknowledgement}
This work has been supported by the German Federal Ministry of Education and Research (BMBF) under
contract number 05M13PMA.

\bibliographystyle{unsrt}
\bibliography{bib}
\end{document}